\documentclass[10pt]{article}
\usepackage{amsmath,amssymb,amsfonts}
\usepackage{epsf}
\usepackage{graphicx}
\usepackage{youngtab}
\usepackage{tikz}
\newtheorem{theorem}{Theorem}[section]
\newtheorem{lemma}{Lemma}[section]
\newtheorem{definition}{Definition}[section]

\newtheorem{corollary}{Corollary}[section]
\newtheorem{remark}{Remark}[section]

\date{}

\title{Large Parts of Random Plane Partitions: a Poisson Limit Theorem}
\bigskip

\author{{\bf Ljuben Mutafchiev} \thanks{This work was partially supported by Project KP-06-N32/8 with the Bulgarian Ministry of
  Education and Science}\\
American University in Bulgaria, 2700 Blagoevgrad, Bulgaria \\ and
Institute of Mathematics and Informatics of the \\ Bulgarian
Academy of Sciences
\\ \tt {ljuben@aubg.edu; tel. +359 2 9435751}}

\begin{document}
\maketitle

\begin{abstract}
 We propose an approach for asymptotic analysis of plane
 partition statistics related to counts of parts whose sizes
 exceed a certain suitably chosen level. In our study, we use the
 concept of conjugate trace of a plane partition of the positive
 integer $n$, introduced by Stanley in 1973. We derive
 generating functions and determine the asymptotic behavior of counts
 of large parts using a general scheme based on the saddle
 point method. In this way, we are able to prove a Poisson limit
 theorem for the number of parts of a random and
 uniformly chosen plane partition of $n$, whose sizes are greater than
 a function $m=m(n)$ as $n\to\infty$. An explicit expression for $m(n)$ is also
 given.
\end{abstract}

\vspace{.5 cm}



 \vspace{.2cm}

\section{Introduction, Motivation and Statement of the Main Result}

Plane partitions were originally introduced by Young \cite{Y01} as
a natural generalization of integer partitions in the plane.
Enumerative problems for plane partitions were first studied via
generating functions by MacMahon \cite{M12} (see also \cite{M16}).
To describe the problem, we will introduce first the concept of a
linear integer partition. For a positive integer $n$, by a
partition $\lambda=(\lambda_1,\lambda_2,\dots,\lambda_k)$ of $n$,
we mean the representation
\begin{equation}\label{ip}
n=\lambda_1+\lambda_2+\dots+\lambda_k,
\end{equation}
where $k\ge 1$ and the integers $\lambda_j, j=1,2,\dots,k$, are
arranged in non-increasing order:
$\lambda_1\ge\lambda_2\ge\dots\ge\lambda_k>0$. The summands
$\lambda_j$ in (\ref{ip}) are usually called parts of $\lambda$.
The Ferrers diagram of a partition is an array of boxes (or cells)
in the plane, left-justified, with $\lambda_j$ boxes in the $j$th
row counting from the bottom. Reading consecutively the numbers of
cells in the columns of the array of the partition $\lambda$,
beginning from the leftmost column, we get the conjugate partition
$\lambda^*=(\lambda_1^*,\lambda_2^*,\dots,\lambda_K^*)$, where
$K=\lambda_1$. For example, the Ferrers diagram of the partition
$\tilde{\lambda}=(5,4,3,3,2,2,2,1)$ of $n=22$ as
$22=5+4+3+3+2+2+2+1$ and its conjugate partition
$\tilde{\lambda}^*=(8,7,4,2,1)$ (that is, $22=8+7+4+2+1$) are
presented in Figure 1.

$$
\centering \yng(1,2,2,2,3,3,4,5)\qquad \yng(1,2,4,7,8)
$$
$$
\text{Figure 1: Ferrers diagrams of $\tilde{\lambda}$ and
$\tilde{\lambda}^*$}
$$

    Let $p(n)$ denote the total number of integer partitions of
    $n\ge 1$. For the generating function
    $$
    P(x)=1+\sum_{n=1}^\infty p(n)x^n
    $$
    of the sequence $\{p(n)\}_{n\ge 1}$, Euler established the
    following identity:
    \begin{equation}\label{pgf}
    P(x)=\prod_{j=1}^\infty (1-x^j)^{-1};
    \end{equation}
see, e.g., \cite[Chapter~1]{A76}. Hardy and Ramanujan \cite{HR18}
developed the so called circle method and applied it to an
asymptotic analysis for the coefficients of $P(x)$. In this way,
they determined asymptotically the numbers $p(n)$ as follows:
$$
p(n)
\sim\frac{1}{4n\sqrt{3}}\exp{\left(\pi\sqrt{\frac{2n}{3}}\right)},
\quad n\to\infty.
$$
For more details and a more precise asymptotic expansion for
$p(n)$, we refer the reader to \cite{R37} and \cite[Chapter
5]{A76}.

The planar analogue of (\ref{ip}) is called a plane partition. A
plane partition $\omega$ of the positive integer $n$ is an array
of non-negative integers
\begin{equation}\label{pp}
\begin{array}{cccccccccccc}
\omega_{1,1} & \omega_{1,2} & \omega_{1,3} \quad \cdots \\
\omega_{2,1} & \omega_{2,2} & \omega_{2,3} \quad \cdots \\
\cdots & \cdots & \cdots \\
\end{array}
\end{equation}
that satisfy $\sum_{h,j\ge 1}\omega_{h,j}=n$, and the rows and
columns in (\ref{pp}) are arranged in non-increasing order:
$\omega_{h,j}\ge\omega_{h+1,j}$ and
$\omega_{h,j}\ge\omega_{h,j+1}$ for all $h,j\ge 1$. The non-zero
entries $\omega_{h,j}>0$ are called parts of $\omega$. If there
are $\lambda_h$ parts in the $h$th row, so that for some $l$,
$\lambda_1\ge\lambda_2\ge\dots\ge\lambda_l >\lambda_{l+1}=0$, then
the (linear) partition
$\lambda=(\lambda_1,\lambda_2,\dots,\lambda_l)$ of the integer
$s=\lambda_1+\lambda_2+\dots+\lambda_l$ is called the shape of
$\omega$. We also say that $\omega$ has $l$ rows and $s$ parts.
Sometimes, for the sake of brevity, the zeroes in array (\ref{pp})
are deleted. For instance, the abbreviation
\begin{equation}\label{ep}
\begin{array}{cccc}
5 & 4 & 1 & 1 \\
3 & 2 & 1 \\
2 & 1 \\
\end{array}
\end{equation}
represents a plane partition $\tilde{\omega}$ of $n=20$ with $l=3$
rows and $s=9$ parts. Any plane partition $\omega$ has an
associated solid diagram $\Delta=\Delta(\omega)$ of volume $n$,
which is considered as three-dimensional analogue of the Ferrers
diagram of a linear integer partition. It is defined as a set of
$n$ integer lattice points
$\mathbf{x}=(x_1,x_2,x_3)\in\mathbb{N}^3$, such that if
$\mathbf{x}\in\Delta$ and $x_j^\prime\le x_j, j=1,2,3$, then
$\mathbf{x}^\prime=(x_1^\prime,x_2^\prime,x_3^\prime)\in\Delta$
too. (Here $\mathbb{N}$ denotes the set of all positive integers.)
Indeed, the entry $\omega_{h,j}$ can be interpreted as the height
of the column of unit cubes stacked along the vertical line
$x_1=h, x_2=j$, and the solid diagram is the union of all such
columns.

Let $q(n)$ denote the total number of plane partitions of the
positive integer $n$ (or, the total number of solid diagrams of
volume $n$). A basic generating function identity established by
MacMahon \cite{M12} implies that the generating function of the
sequence $\{q(n)\}_{n\ge 1}$,
$$
Q(x)=1+\sum_{n=1}^\infty q(n)x^n,
$$
satisfies
\begin{equation}\label{qgf}
Q(x)=\prod_{j=1}^\infty(1-x^j)^{-j}
\end{equation}
(more details may be also found, e.g., in
\cite[Corollary~18.2]{S71} and \cite[Corollary~11.3]{A76}). The
asymptotic form of the numbers $q(n)$, as $n\to\infty$, has been
obtained by Wright \cite{W31} (see also \cite{MK06} for a little
correction). It is given by the following formula:
\begin{equation}\label{wr}
q(n)\sim \frac{(\zeta(3))^{7/36}}{2^{11/36}(3\pi)^{1/2}} n^{-25/36}
\exp{(3(\zeta(3))^{1/3}(n/2)^{2/3}+2\gamma)},
\end{equation}
where
$$
\zeta(z)=\sum_{j=1}^\infty j^{-z}
$$
is the Riemann zeta function and
$$
\gamma=\int_0^\infty\frac{u\log{u}}{e^{2\pi u}-1}du
=\frac{1}{2}\zeta^\prime(-1).
$$
 (The constant $\zeta^\prime(-1)=-0.1654\dots$ is closely related to the
Glaisher-Kinkelin constant; see \cite{F03}).

\begin{remark}. In fact, Wright \cite{W31} proposed a variant
of the circle method, which allows him to obtain an asymptotic
expansion for $q(n)$.
\end{remark}

Further, we need the concepts of conjugate trace and trace of a
plane partition, introduced by Stanley \cite{S73}.

\begin{definition}.\label{d1}
The conjugate trace of the plane partition $\omega$, given by the
array (\ref{pp}), is defined to be the number of parts
$\omega_{h,j}$ of $\omega$ satisfying $\omega_{h,j}\ge h$.
\end{definition}

Hence the conjugate trace of the plane partition in example
(\ref{ep}) is $6$.

\begin{definition}.\label{d2}
The trace of a plane partition $\omega$, given by (\ref{pp}), is
defined to be the sum $\sum_h\omega_{h,h}$.
\end{definition}

Let $T_{lt}^*(n)$ ($T_{lt}(n)$) be the number of plane partitions
of $n$ with $\le l$ rows and conjugate trace (trace) $t$. Stanley
\cite{S73} applied a bijection, established by Bender and Knuth
\cite{BK72}, and showed that
\begin{equation}\label{tr}
T_{lt}^*(n)=T_{lt}(n).
\end{equation}
Setting $T_t^*(n)=\lim_{l\to\infty}T_{lt}^*(n)$ and
$T_t(n)=\lim_{l\to\infty}T_{lt}(n)$, he obtained the following
identities:
\begin{equation}\label{stangf}
1+\sum_{n=1}^\infty\sum_{t=1}^\infty T_t^*(n)y^t x^n
=1+\sum_{n=1}^\infty\sum_{t=1}^\infty T_t(n)y^t x^n
=\prod_{j=1}^\infty (1-yx^j)^{-j}.
\end{equation}
We notice that a linear (one-dimensional) partition $\lambda$ has
$2!=2$ aspects - partition $\lambda$ itself and its conjugate
partition $\lambda^*$, while in the case of plane partitions, we
observe $3!=6$ aspects obtained from the six permutations of the
three axis in the solid diagram. (Here we prefer to use MacMahon's
term "aspect" \cite[Section 427]{M16} rather than "conjugate" used
by Stanley \cite[p. 58]{S73}.) Stanley \cite[Section 3]{S73}
showed that every plane partition $\omega$ of the positive integer
$n$ whose conjugate trace is $t$ has exactly one aspect
$\omega^\prime$  with the same number of rows and trace equal to
$t$. Suppose that the rows of $\omega$ and $\omega^\prime$ are
numbered from top to bottom by $1,2,\dots$. The unique partition
$\omega^\prime$ is obtained from $\omega$ by taking the linear
conjugate partition of row $k$ of $\omega$ and writing then its
parts in a non-increasing order on row $k$ of $\omega^\prime$.
This correspondence explains why both (\ref{tr}) and
(\ref{stangf}) hold. For example, using the conjugates of the
partitions $11=5+4+1+1, 6=3+2+1,$ and $3=2+1$ of the rows in the
partition $\tilde{\omega}$ displayed by (\ref{ep}), we obtain
$\tilde{\omega}^\prime$ as
$$
\begin{array}{ccccc}
4 & 2 & 2 & 2 & 1 \\
3 & 2 & 1 \\
2 & 1, \\
\end{array}
$$
whose trace is obviously $6$ and is equal to the conjugate trace
of $\tilde{\omega}$.

Let $\Omega(n)$ be the set of all plane partitions of $n$, and let
$\Lambda(n)$ be the set of all linear integer partitions of $n$.
We introduce the uniform probability measures $\mathbb{P}$ and
$\mathcal{P}$ on these two sets, respectively. That is, we assign
the probability $1/q(n)$ to each plane partition of $n$ and the
probability $1/p(n)$ to each linear partition of $n$. In this way,
each numerical characteristic of a plane partition from
$\Omega(n)$ and of a linear partition from $\Lambda(n)$ becomes a
random variable (or, a statistic in the sense of the random
generation of plane or linear partitions of $n$). The analysis of
linear integer partitions in terms of probabilistic limit theorems
was initiated by Erd\H{o}s and Lehner \cite{EL41} who found an
appropriate normalization for the largest part (for the number of
parts, by the conjugation of the Ferrers diagram) in a random
partition from $\Lambda(n)$ and established a weak convergence to
the extreme value (Gumbel) distribution as $n\to\infty$.
Subsequent work in this direction was continued by many authors.
Special interest in these studies was to determine the asymptotic
behavior of counts of big part sizes of a random partition from
$\Lambda(n)$ (say, part sizes greater than a certain suitable
value $m=m(n)$). For typical results, we refer the reader to
\cite{ST77}, \cite{F93} and \cite{P97}. In \cite{F93}, among other
important results, Fristedt proved the following Poisson limit
theorem.

\begin{theorem}. \cite[p. 713]{F93}.
\label{t1} Let $\xi_{m,n}$ be the number of parts greater than $m$
in a random integer partition from the set $\Lambda(n)$, equipped
with the uniform probability measure $\mathcal{P}$. Then, with
respect to $\mathcal{P}$, $\xi_{m,n}$ has a limiting Poisson
distribution with expectation $e^{-c}$, as $n\to\infty$ if, for
any $c\in\mathbb{R}$,
\begin{equation}\label{mpf}
m=\frac{\sqrt{6n}}{\pi} \left(\log{\frac{\sqrt{6n}}{\pi}}
+c\right).
\end{equation}
\end{theorem}

The main purpose of this paper is to establish an analogue of
Theorem \ref{t1} for plane partitions from the set $\Omega(n)$,
equipped with the uniform probability measure $\mathbb{P}$. To
introduce a statistic with an asymptotic behavior similar to that
of $\xi_{m,n}$, in the two-dimensional case we need to take into
account the order of the parts of a plane partition in both
directions: from top to bottom (in rows) and from left to right
(in columns). This idea is, in fact, what Stanley \cite{S73}
stated in the definition of a conjugate trace of a plane partition
(see Definition \ref{d1}). The definition of the statistic that we
propose as an analogue of $\xi_{m,n}$ is given below.

For $0\le m<n$, let $X_{m,n}=X_{m,n}(\omega)$ be the number of
parts of $\omega\in\Omega(n)$ that satisfy the inequalities
$\omega_{h,j}\ge h$ and $\omega_{h,j}>m$. In example (\ref{ep}),
we have $X_{0,20}(\tilde{\omega})=6$, $X_{1,20}(\tilde{\omega})=4,
X_{2,20}(\tilde{\omega})=3$ and $X_{3,20}(\tilde{\omega})=2$. Now,
we state our main result.

\begin{theorem}. \label{t2} With respect to the probability measure
$\mathbb{P}$, the random variable $X_{m,n}$ has a limiting Poisson
distribution with expectation $\frac{2}{3}e^{-c}$, as $n\to\infty$
if, for any $c\in\mathbb{R}$,
\begin{equation}\label{mp}
m=\left(\frac{n}{2\zeta(3)}\right)^{1/3}
\left(\log{\left(\frac{n}{2\zeta(3)}\right)^{2/3}} +\log{\log{n}}
+c\right).
\end{equation}
\end{theorem}

\begin{remark}. In \cite{KM07} it is shown that the trace of
a random plane partition of $n$, appropriately normalized,
converges in distribution to a standard normal random variable as
$n\to\infty$. We notice that the conjugate trace and the trace of
a random plane partition of $n$ have one and the same probability
distribution with respect to the probability measure $\mathbb{P}$.
This follows from Stanley's one-to-one correspondence \cite{S73},
described above, and his identity (\ref{stangf}). Hence, the limit
theorem in \cite{KM07} is also valid for the conjugate trace of a
plane partition.
\end{remark}

In the proof of Theorem \ref{t2} we follow a generating function
approach. Let $\mathbb{E}$ denote the expectation taken with
respect to the probability measure $\mathbb{P}$ on the space
$\Omega(n)$. We observe that the generating function of the
expectations $\{\mathbb{E}(y^{X_{m,n}})\}_{n\ge 1}$ satisfies an
identity whose right-hand side is of the form $Q(x)f_m(x,y)$,
where $Q(x)$ is given by (\ref{qgf}) and the function $f_m(x,y)$
will be specified later. Then, we apply the saddle point method in
a form given by Hayman \cite{H56} (see also, e.g., \cite[Chapter
VIII.5]{FS09}).

Finally, we consider the simpler statistic
$Z_{m,n}=Z_{m,n}(\omega)$ counting the number of parts which are
greater than $m$ in a randomly chosen $\omega\in\Omega(n)$. It is
not difficult to show that if $m=m(n)$ is given by (\ref{mp}),
then the difference $Z_{m,n}-X_{m,n}$ tends to $0$ in probability
as $n\to\infty$. Hence Theorem \ref{t2} implies the following
corollary.

\begin{corollary}. \label{c1} With respect to the probability measure
$\mathbb{P}$, the random variable $Z_{m,n}$ has a limiting Poisson
distribution with expectation $\frac{2}{3}e^{-c}$ as $n\to\infty$
if $m$ satisfies (\ref{mp}).
\end{corollary}

Our paper is organized as follows. In Section 2 we include the
necessary generating function identities and the asymptotic
results that will be used further. The proofs of Theorem \ref{t2}
and Corollary \ref{c1} are given in Section 3. Some concluding
remarks are given in Section 4.

\section{Preliminary Results}
\setcounter{equation}{0}

Consider a plane partition $\omega\in\Omega(n)$, defined by the
array (\ref{pp}). Let $L_n=L_n(\omega)$ be the largest part size
of $\omega$ and let $R_n=R_n(\omega)$ be the number of rows in it.
Suppose that, for a certain $\omega$, $R_n\le s$ and define the
subsets of parts of $\omega$ by $\{\omega_{h,j}:\omega_{h,j}=k\ge
h\}$ for $k=1,2,\dots,n$ (possibly some of the last subsets are
empty). Let $Y_{k,n}=|\{\omega_{h,j}:\omega_{h,j}=k\ge h\}|$,
where by $|A|$ we denote the cardinality of the set $A$. The next
lemma extends the result of Theorem 2.2 from \cite{S73}. We state
it in terms of probability generating functions. We assume there
that the randomly chosen $\omega\in\Omega(n)$ is such that
$R_n(\omega)\le s$ and $L_n(\omega)\le l$, for fixed $s$ and $l$
(i.e. the intersection $\{\omega: R_n\le s\}\cap\{\omega: L_n\le
l\}$ is non-empty). For an arbitrary random variable
$U_n=U_n(\omega), \omega\in\Omega(n)$, restricted on $\{\omega:
R_n\le s\}\cap\{\omega: L_n\le l\}$, by $\mathbb{E}(U_n, R_n\le s,
L_n\le l)$ we denote its expectation. (Obviously, after the two
passages to the limit: $s\to\infty$ and $l\to\infty$, we will
obtain the expectation of $U_n$, with respect to the probability
measure $\mathbb{P}$ on the whole $\Omega(n)$, that is
$\mathbb{E}(U_n)$.)

\begin{lemma}. \label{l1} We have
\begin{eqnarray}\label{gf}
& & 1+\sum_{n=1}^\infty q(n)\mathbb{E}
(y_1^{Y_{1,n}}y_2^{Y_{2,n}}\dots y_n^{Y_{n,n}}, R_n\le s, L_n\le
l)
x^n \nonumber \\
& & =\prod_{k=1}^s\prod_{j=1}^l (1-y_j x^{k+j-1})^{-1},
\end{eqnarray}
where $x,y_1,\dots,y_n$ are formal variables.
\end{lemma}

{\it Sketch of the proof.} We notice first that Definition
\ref{d1} implies that $\sum_{j\ge 1}Y_{j,n}$ equals the conjugate
trace of a plane partition. Stanley \cite[Theorem 2.2]{S73} showed
that from (\ref{tr}) it follows that
\begin{eqnarray}\label{sgf}
& & 1+\sum_{n=1}^\infty q(n) \mathbb{E}(y^{\sum_{j\ge 1}Y_{j,n}},
R_n\le s, L_n\le l)x^n \nonumber \\
& & =\prod_{k=1}^s\prod_{j=1}^l (1-yx^{k+j-1})^{-1},
\end{eqnarray}
where $x$ and $y$ are formal variables. Clearly, (\ref{gf}) is a
slight extension of (\ref{sgf}), in which the contribution of the
parts $k$ are separated by the variables $y_k$. The proof of
(\ref{gf}) follows the same line of reasoning as in \cite{S73}. It
is based on two bijections. The first one, due to Knuth
\cite{K70}, is as follows:

(K) There is a one-to-one correspondence between ordered pairs
$(\omega_1,\omega_2)$ of column strict plane partitions of the
same shape and matrices $(b_{jk})_{j,k\ge 1}$ of non-negative
integers. In this correspondence, (i) the number $k$ appears in
$\omega_1$ exactly $\sum_j b_{jk}$ times, and (ii) $k$ appears in
$\omega_2$ exactly $\sum_j b_{kj}$ times.

(By a column strict plane partition we mean a plane partition
whose non-zero entries are strictly decreasing in each column.)

The second bijection is given by Bender and Knuth \cite{BK72}.

(BK) There is a one-to-one correspondence between the plane
partition $\omega\in\Omega(n)$ with $R_n(\omega)=s, L_n(\omega)=l$
and $Y_{j,n}(\omega)=t_j, j=1,2,\dots,n$, with $t=\sum_{j=1}^n
t_j$, and pairs $(\omega_1,\omega_2)$ of column strict partitions,
so that the largest part of $\omega_2$ is $s$, the largest part of
$\omega_1$ is $l$, the number of parts in the $j$th row of
$\omega_1$ or $\omega_2$ is $t_j$ and the conjugate trace $t$ of
$\omega$ equals the total number of parts $\sum_{j=1}^n t_j$ of
$\omega_1$ or $\omega_2$. Moreover, if $\omega_k$ is a partition
of $n_k, k=1,2$, then $n=n_1+n_2-t$; see also \cite[p. 57]{S73}.

In this way, using bijection (BK), we establish that the count
$q(n)\mathbb{P}(Y_{j,n}=t_j,j=1,2,\dots,n, R_n\le s,L_n\le l)$ is
equal to the number of pairs $(\omega_1,\omega_2)$ of column
strict plane partitions of the same shape satisfying: (i) the
largest part of $\omega_1$ is $\le s$, (ii) the largest part of
$\omega_2$ is $\le l$, (iii) the number of parts in row $j$ of
$\omega_1$ or $\omega_2$ is $t_j$, (iv) the number of parts of
$\omega_1$ or $\omega_2$ is $t=\sum_{j=1}^n t_j$, (v) the sum of
the parts of $\omega_1$ and $\omega_2$ is $n+t$. Then, one can
obtain (\ref{gf}), using bijection (K) and following the same
argument as in \cite{S73}. $\Box$

\begin{remark}. We notice that, in a similar way as Stanley did in
\cite{S73}, one can obtain from (\ref{gf}) as corollaries the
following identities after both passages to the limit:
$s\to\infty$ and $l\to\infty$. For instance, we have
\begin{equation}\label{gff}
1+\sum_{n=1}^\infty q(n)\mathbb{E}(\prod_{k=1}^\infty
y_k^{Y_{k,n}}) x^n =\prod_{j=1}^\infty (1-y_j x^j)^{-j},
\end{equation}
which generalizes Stanley's formula (6) in \cite[p. 59]{S73}.
Clearly, (\ref{stangf}) follows from (\ref{gff}) after the
substitution $y_j=y, j=1,2,\dots$.
\end{remark}

Now, recall that $X_{m,n}=\sum_{j>m}Y_{j,n}$. Setting in
(\ref{gff}) $y_1=\dots=y_{[m]}=1$ and $y_j=y$ for $j>m$, where
$[m]$ denotes the integer part of $m$, we obtain
\begin{eqnarray}\label{xgf}
& & 1+\sum_{n=1}^\infty q(n)\mathbb{E}(y^{X_{m,n}}) =\prod_{j\le
m}(1-x^j)^{-j}\prod_{j>m}(1-yx^j)^{-j} \nonumber \\
& & =Q(x)f_m(x,y).
\end{eqnarray}
Here $Q(x)$ is defined by (\ref{qgf}) and
\begin{equation}\label{f}
f_m(x,y)=\prod_{j>m}\left(\frac{1-x^j}{1-yx^j}\right)^j.
\end{equation}

We are now ready to proceed with the preliminaries of our further
asymptotic analysis. We have to study the behavior of the
coefficient $[x^n]Q(x)f_m(x,y)$ of $x^n$ in the Taylor expansion
of the product $Q(x)f_m(x,y)$ as $n\to\infty$, where $Q(x)$ and
$f_m(x,y)$ are defined by (\ref{qgf}) and (\ref{f}), respectively.
We express this coefficient by means of Cauchy's integral formula
using a suitably chosen closed curve around $0$ as a contour of
integration. Since the unit circle is a natural boundary for
$Q(x)$, this contour lies inside the unit disk. We will estimate
the Cauchy integral using a general theorem due to Hayman
\cite{H56} whose proof is based on the saddle point method. We
will describe next the wide class of functions to which Hayman's
theorem apply. We employ the terminology given in \cite[Chapter
VIII.5]{FS09}.

Consider a function $G(x)=\sum_{n=0}^\infty g_n x^n$ that is
analytic for $|x|<\rho, 0<\rho\le\infty$. For $0<r<\rho$, we set
\begin{equation}\label{a}
a(r)=r\frac{G^\prime(r)}{G(r)},
\end{equation}
\begin{equation}\label{b}
b(r)=\frac{rG^\prime(r)}{G(r)}+r^2\frac{G^{\prime\prime}(r)}{G(r)}
-r^2\left(\frac{G^\prime(r)}{G(r)}\right)^2.
\end{equation}
We assume that $G(x)>0$ for $x\in (\rho_0,\rho)\subset (0,\rho)$
and satisfies the following three conditions:

{\it Capture condition.} $lim_{r\to\rho} a(r)=\infty$ and
$\lim_{r\to\rho} b(r)=\infty$.

{\it Locality condition.} For some function $\delta=\delta(r)$
defined over $(\rho_0,\rho)$ and satisfying $0<\delta<\pi$, one
has
$$
G(r e^{i\theta})\sim G(r) e^{i\theta a(r)-\theta^2 b(r)/2}
$$
 as $r\to\rho$, uniformly for $|\theta|\le\delta(r)$.

 {\it Decay condition.}
$$
 G(r e^{i\theta})
=o\left(\frac{G(r)}{\sqrt{b(r)}}\right)
$$
 as $r\to\rho$ uniformly for $\delta(r)<|\theta|\le\pi$.

 \begin{definition}. \label{dhayman} A function $G(x)$ which satisfies the capture,
 locality and decay conditions is called admissible in the sense
 of Hayman.
 \end{definition}

 {\bf Hayman's Theorem.} {\it Let $G(x)$ be a Hayman admissible function
and $r=r_n$ be the unique solution of the equation
\begin{equation}\label{arn}
a(r)=n.
\end{equation}
Then the Taylor coefficients $g_n$ of $G(x)$ satisfy, as
$n\to\infty$,
\begin{equation}\label{gn}
g_n\sim\frac{G(r_n)}{r_n^n\sqrt{2\pi b(r_n)}}
\end{equation}
with $b(r_n)$ given by (\ref{b}).}

Hayman's theorem was applied in \cite{M18} to obtain a general
asymptotic estimate for $[x^n]Q(x)F(x)$, where $Q(x)$ is defined
by (\ref{qgf}) and $F(x)$ is suitably restricted on its behavior
in any neighborhood of $x=1$. This result was then used to derive
asymptotics of the expectations of several plane partition
statistics. For linear integer partition statistics, a method for
the asymptotic analysis of $[x^n]P(x)F(x)$, where $P(x)$ is
defined by (\ref{pgf}), was developed by Grabner et al.
\cite{GKW14}. In our asymptotic analysis of $[x^n]Q(x)f_m(x,y)$ we
are not able to apply directly the general result from
\cite[Theorem 1]{M18}. In fact, the function $f_m(x,y)$ (see
(\ref{f})) is analytic only for $|x|<1$, while $F(x)$ in
\cite{M18} satisfies a more general assumption in a neighborhood
of the point $x=1$. In addition, $f_m(x,y)$ depends on the growth
of $n$ since the parameter $m=m(n)$ satisfies relation (\ref{mp}).
Finally, the dependency on a second variable $y$ in $f_m(x,y)$
plays an important role in our study since the estimates that we
will explore further have to be uniform for $y\in [0,1)$. Hence
our further asymptotic analysis relies on the observations for the
MacMahon's generating function $Q(x)$ established in \cite{M18}.
The asymptotic behavior of $f_m(x,y)$, as $x\to 1$ and $|x|<1$, is
studied separately at the end of this section.

We continue with a lemma related to $Q(x)$ (see (\ref{qgf})). For
its proof, we refer the reader to \cite[pp. 261-265]{M18}.
Further on, in (\ref{a}), (\ref{b}) and (\ref{gn}) we set
$G(x):=Q(x)$, $\rho:=1$ and $g_n:=q(n)$.

\begin{lemma}. \label{l2} (i) The unique solution
of (\ref{arn}) is $r=r_n=e^{-d_n}$, where, for large $n$, the
sequence $\{d_n\}_{n\ge 1}$ has the following expansion:
\begin{equation}\label{d}
d_n=\left(\frac{2\zeta(3)}{n}\right)^{1/3} -\frac{1}{36n}
+O(n^{-1-\beta})
\end{equation}
and $\beta>0$ is a certain fixed constant. Moreover, as
$n\to\infty$,
\begin{equation}\label{bd}
b(e^{-d_n})\sim\frac{3n^{4/3}}{2\zeta^{1/3}(3)},
\end{equation}
where $b(r)$ is defined by (\ref{b}). Hence (\ref{arn}) and
(\ref{bd}) imply that $Q(x)$ satisfies Hayman's "capture"
condition as $n\to\infty$.

(ii) With $d_n$ and $b(e^{-d_n})$ as in part (i), we have
\begin{equation}\label{locq}
e^{-i\theta n}\frac{Q(e^{-d_n+i\theta})}{Q(e^{-d_n})}
=e^{-\theta^2 b(e^{-d_n})/2}(1+O(1/\log^3{n}))
\end{equation}
as $n\to\infty$ uniformly for $|\theta|\le\delta_n$, where
\begin{equation}\label{delta}
\delta_n=\frac{d_n^{5/3}}{\log{n}} =\frac{1}{\log{n}}
\left(\frac{2\zeta(3)}{n}\right)^{5/9} (1+O(n^{-2/3})).
\end{equation}
(The last equality follows from (\ref{d}).) In addition,
(\ref{locq}) and (\ref{delta}) show that Hayman's "locality"
condition holds for $Q(x)$ with $\delta_n:=\delta(e^{-d_n})$.

(iii) For sufficiently large $n$, we have
\begin{eqnarray}\label{decq}
& & |Q(e^{-d_n+i\theta})|\le Q(e^{-d_n})e^{-Cd_n^{-2/3}} \nonumber
\\
& & \le Q(e^{-d_n}) e^{-C^\prime n^{2/9}/\log^2{n}}
=o\left(\frac{Q(e^{-d_n})}{\sqrt{b(e^{-d_n})}}\right)
\end{eqnarray}
uniformly for $\delta_n\le |\theta|\le\pi$, where $C, C^\prime>0$
are absolute constants and $d_n$ and $\delta_n$ satisfy (\ref{d})
and (\ref{delta}), respectively. By (\ref{decq}) $Q(x)$ satisfies
also Hayman's "decay" condition with $\delta_n$ as in part (ii).
\end{lemma}

\begin{remark}. Clearly, for sufficiently large $n$, the arc
$(-\delta_n,\delta_n)$ on the circle $x=e^{-d_n+i\theta},
-\pi<\theta\le\pi,$ becomes close to the main singularity $x=1$ of
MacMahon's generating function $Q(x)$. Furthermore, (\ref{locq})
and (\ref{decq}) show that $Q(e^{-d_n+i\theta})$ significantly
changes its behavior when $\theta$ leaves the interval
$(-\delta_n,\delta_n)$. Moreover, Lemma \ref{l2} and Definition
\ref{dhayman} show that MacMahon's generating function $Q(x)$ is
admissible in the sense of Hayman. Therefore, by Hayman's theorem
we have
\begin{equation}\label{wh}
q(n)\sim\frac{e^{nd_n}Q(e^{-d_n})}{\sqrt{2\pi b(e^{-d_n})}}, \quad
n\to\infty,
\end{equation}
where $d_n$ and $b(e^{-d_n})$ are given by (\ref{d}) and
(\ref{bd}), respectively. In the Appendix of \cite{M18} it is
shown that (\ref{wh}) implies the corrected form of Wright's
formula (\ref{wr}).
\end{remark}

Our last task in this section is to study the behavior of
$f_m(e^{-id_n+i\theta},y)$; see (\ref{f}). We obtain uniform
estimates for the following two cases: $\theta$ belongs to a
neighborhood of $0$, and $\theta$ varies arbitrarily in the whole
interval $(-\pi,\pi]$.

\begin{lemma} \label{l3} (i) If $d_n$ and $m=m(n)$ satisfy (\ref{d}) and
(\ref{mp}), respectively, then
\begin{equation}\label{limf}
\lim_{n\to\infty}\frac{f_m(e^{-id_n+i\theta},y)}{f_m(e^{-d_n},y)}=1
\end{equation}
uniformly for $|\theta|\le\delta_n$ and $y\in [0,1)$, where
$\delta_n$ is given by (\ref{delta}).

(ii) Let $d_n$ and $m$ be the same as in part (i). Then, for any
$\theta\in (-\pi,\pi]$ and sufficiently large $n$, we have
\begin{equation}\label{fo}
f_m(e^{-id_n+i\theta},y)=O(1)
\end{equation}
uniformly for $y\in [0,1)$.
\end{lemma}

{\it Proof.} (i) We let $\log{x}$ denote the main branch of the
logarithmic function, that is, we assume that $\log{x}<0$ if
$0<x<1$. Next, using (\ref{f}), we represent the function
$f_m(x,y)$ in the following way:
\begin{equation}\label{fone}
f_m(x,y) =\exp{(\sum_{j>m}jg_j(x,y))},
\end{equation}
where
$$
g_j(x,y)=\log{\frac{1-x^j}{1-yx^j}}.
$$
By the Taylor formula with $x=e^{-d_n+i\theta}$, we have
\begin{equation}\label{gj}
g_j(e^{-d_n+i\theta},y)=g_j(e^{-d_n},y) +O\left(\delta_n|
\frac{\partial}{\partial x} g_j(x,y)|_{x=e^{-d_n}}\right),
\end{equation}
since, for any $\theta_0\in (-\delta_n,\delta_n)$, we have
$|e^{i\theta_0} -1|\le |\theta_0|<\delta_n$. A simple calculation
shows that
\begin{equation}\label{gjder}
\frac{\partial}{\partial x} g_j(x,y) =\frac{jx^{j-1} (y-1)}
{(1-x^j)(1-yx^j)}.
\end{equation}
Hence, using (\ref{fone}) - (\ref{gjder}), we can write
\begin{equation}\label{esmn}
\frac{f(e^{-d_n+i\theta},y)}{f(e^{-d_n},y)} =e^{S_{m,n}},
\end{equation}
where the sum $S_{m,n}$ of the remainder terms in the Taylor
expansions of $g_j$ satisfies the following estimate:
\begin{equation}\label{smn}
S_{m,n}=O\left(\delta_n(1-y) \sum_{j>m} \frac{j^2
e^{-jd_n}}{(1-e^{-jd_n})(1-y e^{-jd_n})}\right)
\end{equation}
uniformly for $|\theta|<\delta_n$ and $y\in [0,1)$. The sum on the
right-hand side of (\ref{smn}) can be interpreted as a Riemann sum
with step size $d_n$. At this moment, it is more convenient for us
to express $m$ as a function of $d_n$, given by (\ref{d}). We set
\begin{equation}\label{mpd}
m=d_n^{-1}(\log{d_n^{-2}} +\log{\log{n}} +c), \quad
c\in\mathbb{R}.
\end{equation}
By (\ref{mpd}) the lower bound of the integral is
$md_n=\log{d_n^{-2}}+\log{\log{n}}+c$. Thus we have
\begin{eqnarray}\label{summ}
& & \sum_{j>m} \frac{j^2 e^{-jd_n}}{(1-e^{-jd_n})(1-y e^{-jd_n})}
\nonumber \\
& &=O\left(d_n^{-3}\int_{\log{d_n^{-2}}+\log{\log{n}}+c}^\infty
\frac{u^2 e^{-u}}{(1-e^{-u})(1-ye^{-u})}du\right) \nonumber \\
& &=O\left(\frac{d_n^{-3}}{1-y}
\int_{\log{d_n^{-2}}+\log{\log{n}}+c}^\infty \frac{u^2 e^{-u}}
{1-e^{-u}} du\right).
\end{eqnarray}
The last integral is related to the Debye function of order $2$.
It is easy to see, using formula 27.1.2 in \cite{AS65}, that
\begin{equation}\label{deb}
\int_t^\infty\frac{u^2}{e^u-1}du =(t^2+2t+2)e^{-t} +O(t^2
e^{-2t}), \quad t\to\infty.
\end{equation}
We will also need asymptotic expansions for $d_n^{-1}$, $d_n^{-2}$
and $\log{d_n^{-2}}$. From (\ref{d}) it follows that
\begin{equation}\label{dminus}
d_n^{-1} =\left(\frac{n}{2\zeta(3)}\right)^{1/3}
+\frac{1}{36(2\zeta(3))^{2/3} n^{1/3}} +O(n^{-1/3-\beta}),
\end{equation}

$$
d_n^{-2} =\left(\frac{n}{2\zeta(3)}\right)^{2/3}
+\frac{1}{36\zeta(3)} +O(n^{-\beta})
$$
and
\begin{equation}\label{logd}
\log{d_n^{-2}} =\frac{2}{3}\log{n} -\frac{2}{3}\log{(2\zeta(3))}
+O(n^{-2/3}).
\end{equation}
We notice here that from (\ref{dminus}) and (\ref{logd}) it
follows that
\begin{eqnarray}
& & d_n^{-1}(\log{d_n^{-2}}+\log{\log{n}}+c) \nonumber \\
& & =\left(\frac{n}{2\zeta(3)}\right)^{1/3}
\left(\log{\left(\frac{n}{2\zeta(3)}\right)^{2/3}} +\log{\log{n}}
+c\right)+o(1), \nonumber
\end{eqnarray}
i.e., the difference between the right-hand sides of (\ref{mp})
and (\ref{mpd}) tends to $0$ as $n\to\infty$, which justifies the
value of $m$, given by (\ref{mp}). Moreover, from (\ref{deb}) and
(\ref{logd}) we obtain
\begin{equation}\label{dint}
d_n^{-3}\int_{\log{d_n^{-2}}+\log{\log{n}}+c}^\infty
\frac{u^2}{e^u-1} du =O(d_n^{-1}\log{n}).
\end{equation}
Combining (\ref{d}), (\ref{delta}), (\ref{smn}) and (\ref{dint}),
we see that
$$
S_{m,n}=O(d_n^{2/3})=O(n^{-2/9}),
$$
which implies that the ratio in (\ref{esmn}) tends to $1$, as
$n\to\infty$, uniformly for $|\theta|\le\delta_n$ and $y\in
[0,1)$. This completes the proof of part (i).

(ii) For any $\theta\in(-\pi,\pi]$, we observe that
\begin{eqnarray}\label{fupperb}
& & |f_m(e^{-id_n+i\theta},y)| =\exp{\left(\sum_{j>m}j
\log{\frac{|1-e^{-jd_n\theta}|}{|1-ye^{-jd_n\theta}|}}\right)}
\nonumber \\
& & =\exp{\left(\frac{1}{2}\sum_{j>m}
j\log{\left(\frac{1-2e^{-jd_n}\cos{(j\theta)}+e^{-2jd_n}}
{1-2ye^{-jd_n}\cos{(j\theta)}+y^2 e^{-2jd_n}}\right)}\right)}
\nonumber \\
& & \le\exp{\left(\sum_{j>m}
j\log{\left(\frac{1+e^{-jd_n}}{1-e^{-jd_n}}\right)}\right)},
\end{eqnarray}
where in the inequality on the last line we used that $y\in [0,1)$
and that $-1\le\cos{(j\theta)}\le 1$. Expanding the logarithm on
the right hand side of (\ref{fupperb}) into powers of $e^{-jd_n}$
with $j>m$ and using (\ref{mpd}) and (\ref{logd}), we obtain
\begin{eqnarray}
& & \log{\left(\frac{1+e^{-jd_n}}{1-e^{-jd_n}}\right)}=2e^{-jd_n}
+2\sum_{k\ge 2}\frac{e^{-j(2k-1)d_n}}{2k-1}
\nonumber \\
& & =2e^{-jd_n} +O(e^{-3md_n}) =2e^{-jd_n}(1+O(e^{-2md_n}))
\nonumber \\
& & =2e^{-jd_n}(1+e^{-(\frac{4}{3}\log{n}+2\log{\log{n}}+O(1))})
=2e^{-jd_n}(1+O(n^{-4/3}/\log^2{n})). \nonumber
\end{eqnarray}
Substituting this expression into the right hand side of
(\ref{fupperb}) and approximating again the underlying sum by a
Riemann integral, from (\ref{dminus}) and (\ref{logd}) we get the
estimate:
\begin{eqnarray}\label{finb}
& & |f_m(e^{-id_n+i\theta},y)| \le (1+O(n^{-4/3}\log^2{n}))
\nonumber \\
& & \times\exp{\left(2\left(\left(\frac{n}{2\zeta(3)}\right)^{1/3}
+O\left(\frac{1}{n}\right)\right)\sum_{j>m}(jd_n)e^{-jd_n}\right)} \nonumber \\
& & \le \exp{\left(C_0
n^{1/3}\int_{\frac{2}{3}\log{n}+\log{\log{n}}}^\infty ue^{-u}
du\right)}
\end{eqnarray}
uniformly for $\theta\in (-\pi,\pi]$ and $y\in [0,1)$, where
$C_0>0$ is an absolute constant. Using the asymptotic formula for
the incomplete gamma function \cite[formula 6.5.32]{AS65}, we
obtain
$$
\int_{\frac{2}{3}\log{n}+\log{\log{n}}}^\infty ue^{-u} du
=O(n^{-2/3}).
$$
Replacing this estimate into the right hand side of (\ref{finb}),
we complete the proof of part (ii). $\Box$

\section{Proof of the Main Result}
\setcounter{equation}{0}

\subsection{Proof of Theorem \ref{t2}}

First, we recall the expressions for $d_n$ and $\delta_n$ given by
(\ref{d}) and (\ref{delta}), respectively. Next we apply the
Cauchy coefficient formula to (\ref{xgf}), using the circle
$x=e^{-d_n+i\theta}, -\pi<\theta\le\pi$, as a contour of
integration. Thus we obtain
\begin{eqnarray}\label{intmain}
& & [x^n]Q(x)f_m(x,y) =\frac{e^{nd_n}}{2\pi}\int_{-\pi}^\pi
 Q(e^{-d_n+i\theta}) f_m(e^{-d_n+i\theta},y) e^{-i\theta n}d\theta
 \nonumber \\
& & =J_{1,n}+J_{2,n},
 \end{eqnarray}
 where
 \begin{equation}\label{intone}
 J_{1,n} =\frac{e^{nd_n}}{2\pi}\int_{-\delta_n}^{\delta_n}
Q(e^{-d_n+i\theta}) f_m(e^{-d_n+i\theta}) d\theta,
\end{equation}

\begin{equation}\label{inttwo}
J_{2,n}=\frac{e^{nd_n}}{2\pi}\int_{\delta_n<\mid\theta\mid\le\pi}
Q(e^{-d_n+i\theta}) f_m(e^{-d_n+i\theta},y) d\theta.
\end{equation}

The estimate of $J_{1,n}$ follows from parts (i) and (ii) of Lemma
\ref{l2} and Lemma \ref{l3}. First, in (\ref{intone}) we perform
the following computation:
\begin{eqnarray}\label{estintone}
& & J_{1,n} =\frac{e^{nd_n}Q(e^{-d_n})f_m(e^{-d_n},y)}{2\pi}
\nonumber \\
& & \times\int_{-\delta_n}^{\delta_n}
\left(\frac{Q(e^{-d_n+i\theta})} {Q(e^{-d_n})}\right)
\left(\frac{f_m(e^{-d_n+i\theta},y)}
{f_m(e^{-d_n},y)}\right)e^{-i\theta n} d\theta \nonumber \\
& & =\frac{e^{nd_n}Q(e^{-d_n})f_m(e^{-d_n},y)}{2\pi} \nonumber \\
& & \times\int_{-\delta_n}^{\delta_n} e^{-\theta^2 b(e^{-d_n})/2}
\left(1+O\left(\frac{1}{\log^3{n}}\right)\right)(1+o(1)) d\theta
\nonumber \\
& & \sim\frac{e^{nd_n}Q(e^{-d_n})f_m(e^{-d_n},y)}{2\pi}
\int_{-\delta_n}^{\delta_n} e^{-\theta^2 b(e^{-d_n})/2}d\theta.
\end{eqnarray}
Note that in the second equality we applied (\ref{locq}) (that is,
Hayman's "locality" condition) and (\ref{limf}) (i.e., Lemma
\ref{l3}(i)). Next, in the last integral of (\ref{estintone}) we
substitute $\theta=u/\sqrt{b(e^{-d_n})}$. We observe that
\begin{eqnarray}
& & \int_{-\delta_n}^{\delta_n} e^{-\theta^2 b(e^{-d_n})/2}d\theta
\sim\frac{1}{\sqrt{b(e^{-d_n})}}
\int_{-\delta_n\sqrt{b(e^{-d_n})}}^{\delta_n\sqrt{b(e^{-d_n})}}
e^{-u^2/2} du \nonumber \\
& & \sim\frac{1}{\sqrt{b(e^{-d_n})}}\int_{-\infty}^\infty
e^{-u^2/2} du =\sqrt{\frac{2\pi}{b(e^{-d_n})}}, \quad n\to\infty,
\nonumber
\end{eqnarray}
since by (\ref{bd}) and (\ref{delta})
$$
\delta_n\sqrt{b(e^{-d_n})} \sim\sqrt{3}(2\zeta(3))^{7/18}
\frac{n^{1/9}}{\log{n}}\to\infty, \quad n\to\infty.
$$
Inserting the estimate of the last integral into
(\ref{estintone}), by Wright's formula (\ref{wh}) we obtain
\begin{eqnarray}\label{festintone}
& & J_{1,n} =\frac{e^{nd_n}Q(e^{-d_n})f_m(e^{-d_n},y)} {\sqrt{2\pi
b(e^{-d_n)}}}  \\
& & =q(n)f_m(e^{-d_n},y)+o(f_m(e^{-d_n},y)q(n))
=q(n)f_m(e^{-d_n},y)+o(q(n)), \nonumber
\end{eqnarray}
where in the last equality we have also used (\ref{fo}).

To estimate $J_{2,n}$, we use Lemma \ref{l2}(iii) and Lemma
\ref{l3}(ii). We apply first the inequality given in (\ref{decq})
and combine it with (\ref{fo}). Thus we observe that
\begin{equation}\label{qfb}
|Q(e^{-d_n+i\theta})f_m(e^{-d_n+i\theta})|\le
C^{\prime\prime}Q(e^{-d_n})e^{-C^\prime n^{2/9}/\log^2{n}}
\end{equation}
uniformly for $\delta_n\le|\theta|\le\pi$, where
$C^{\prime\prime}>0$ is a certain constant. From (\ref{inttwo}),
(\ref{qfb}), (\ref{bd}) and (\ref{wh}) it follows that
\begin{eqnarray}\label{estinttwo}
& & |J_{2,n}|\le\frac{e^{nd_n}}{2\pi} \int_{\delta_n
\le|\theta|<\pi} |Q(e^{-d_n+i\theta})f_m(e^{-d_n+i\theta})|
d\theta
\nonumber \\
& & \le\frac{C^{\prime\prime}e^{nd_n}}{\pi} Q(e^{-d_n})
e^{-C^\prime
n^{2/9}/\log^2{n}}(\pi-\delta_n) \nonumber \\
& & =O\left(\frac{e^{nd_n} Q(e^{-d_n})}{\sqrt{2\pi b(e^{-d_n})}}
n^{2/3} e^{-C^\prime n^{2/9}/\log^2{n}}\right)
\nonumber \\
& & =O(q(n)n^{2/3} e^{-C^\prime n^{2/9}/\log^2{n}}).
\end{eqnarray}
Substituting (\ref{festintone}) and (\ref{estinttwo}) into
(\ref{intmain}) - (\ref{inttwo}), we obtain
$$
\frac{1}{q(n)}[x^n]Q(x)f_m(x,y) =f_m(e^{-d_n},y)+o(1)
$$
uniformly for $y\in [0,1)$. Finally, going back to (\ref{xgf}), we
conclude that
\begin{equation}\label{pgfx}
\mathbb{E}(y^{X_{m,n}})=f_m(e^{-d_n},y)+o(1).
\end{equation}

So, to complete the proof of the theorem it remains to study the
asymptotic behavior of $f_m(e^{-d_n},y)$, with $m$ given by
(\ref{mp}).

We will use an alternative representation for $f_m(x,y)$, which
follows from (\ref{f}). We have
\begin{eqnarray}\label{fmas}
& & f_m(x,y) =\exp{(\sum_{j>m} j(\log{(1-x^j)}-\log{(1-yx^j)}))}
\nonumber \\
& & =\exp{((y-1)(\sum_{j>m} jx^j) +K_m(x,y))},
\end{eqnarray}
where
\begin{equation}\label{rm}
K_m(x,y)=\sum_{j>m}\sum_{k\ge 2}\frac{j}{k} ((yx^j)^k-x^{jk}).
\end{equation}
The sum in the exponent of the right-hand side of (\ref{fmas}) can
be estimated using a Riemann sum approximation as in (\ref{finb}).
Here we need a more precise estimate. Setting $x=e^{-d_n}$, we
obtain
\begin{eqnarray}\label{sumas}
& & \sum_{j>m} je^{-jd_n}\sim d_n^{-2}\int_{md_n}^\infty ue^{-u}
du =d_n^{-2}\int_{\log{d_n^{-2}}+\log{\log{n}}+c}^\infty ue^{-u}
du \nonumber \\
& & =-d_n^{-2} ue^{-u}|_{\log{d_n^{-2}}+\log{\log{n}}+c}^\infty
+d_n^{-2}\int_{\log{d_n^{-2}}+\log{\log{n}}+c}^\infty  e^{-u} du
\nonumber \\
& & =(\log{d_n^{-2}}+\log{\log{n}}+c)\frac{e^{-c}}{\log{n}}
+\frac{e^{-c}}{\log{n}} \nonumber \\
& & =\left(\frac{2}{3}\log{n}+\log{\log{n}} +O(1)\right)
\frac{e^{-c}}{\log{n}} +O\left(\frac{1}{\log{n}}\right) \nonumber
\\
& & =\frac{2}{3} e^{-c}
+O\left(\frac{\log{\log{n}}}{\log{n}}\right),
\end{eqnarray}
where in the fourth equality we have also used (\ref{logd}) and
(\ref{mpd}).

Finally, it remains to study the asymptotic behavior of the
remainder term $K_m(e^{-d_n},y)$ in (\ref{fmas}). First, we change
the order of summation in (\ref{rm}), and then we perform some
algebraic computations in order to see that
\begin{eqnarray}\label{rmn}
& & K_m(e^{-d_n},y) =\sum_{k\ge 2}\frac{(y^k-1)e^{-kd_n}}{k}
\left(\frac{me^{-mkd_n}}{1-e^{-kd_n}}
+\frac{e^{-mkd_n}}{(1-e^{-kd_n})^2}\right) \nonumber \\
& & =K_{m,n}^{(1)}+K_{m,n}^{(2)},
\end{eqnarray}
where
\begin{equation}\label{rmnone}
K_{m,n}^{(1)} =\sum_{k\ge 2} \frac{(y^k-1)e^{-kd_n} me^{-mkd_n}}
{k(1-e^{-kd_n})},
\end{equation}

$$
K_{m,n}^{(2)} =\sum_{k\ge 2}\frac{(y^k-1)e^{-kd_n} e^{-mkd_n}}
{k(1-e^{-kd_n})^2}.
$$
From (\ref{mp}) it follows that
$$
e^{-mkd_n} =\left(\frac{d_n^2 e^{-c}}{\log{n}}\right)^k,
$$

$$
1-e^{-kd_n} =kd_n+O(k^2 d_n^2).
$$
Replacing these two equalities into (\ref{rmnone}) and using
(\ref{d}), we conclude that
\begin{equation}\label{rmne}
K_{m,n}^{(1)} =O\left(\sum_{k\ge 2}\frac{d_n^{2k-2}}{k^2}\right)
=O(d_n^2)=O(n^{-2/3})
\end{equation}
uniformly for $y\in [0,1)$. In the same way, we obtain the
estimate for $K_{m,n}^{(2)}$. We have
\begin{equation}\label{rmnf}
K_{m,n}^{(2)} =O\left(\sum_{k\ge
2}\frac{d_n^{2k-2}}{k^3\log^k{n}}\right)
=O\left(\frac{d_n^2}{\log{n}}\right).
\end{equation}
Combining (\ref{rmn}), (\ref{rmne}), (\ref{rmnf}) and (\ref{d}),
we get
\begin{equation}\label{rmnas}
K(e^{-d_n},y)=O(d_n^2)=O(n^{-2/3}).
\end{equation}
Thus (\ref{fmas}), (\ref{sumas}) and (\ref{rmnas}) imply that
$f_m(e^{-d_n},y)$ approaches $\exp{(\frac{2}{3}e^{-c}(y-1))}$ as
$n\to\infty$ uniformly for $y\in [0,1)$. Now, Theorem \ref{t2}
follows from the continuity theorem for probability generating
functions and (\ref{pgfx}), since $\exp{(\frac{2}{3}e^{-c}(y-1))}$
is the generating function of a Poisson distribution with
expectation $\frac{2}{3}e^{-c}$. $\Box$

\subsection{Proof of Corollary \ref{c1}}

 We recall that the variables $X_{m,n}$
and $Z_{m,n}$ are defined on the set $\Omega(n)$ of all plane
partitions of $n$ equipped with the uniform probability measure
$\mathbb{P}$. Let $\omega=(\omega_{h,j})\in\Omega(n)$ be a plane
partition defined by array (\ref{pp}). Since the restriction
$\omega_{h,j}\ge h$ is removed in the definition of
$Z_{m,n}=Z_{m,n}(\omega)$, then, for every $\omega\in\Omega(n)$,
we have $Z_{m,n}(\omega)\ge X_{m,n}(\omega)$. Furthermore, we
consider the following two events:
$A_{m,n}=\{\omega=(\omega_{h,j})\in\Omega(n):
Z_{m,n}(\omega)-X_{m,n}(\omega)>\epsilon\}$, with $\epsilon>0$,
and $B_{m,n}=\{\omega=(\omega_{h,j})\in\Omega(n):$ there is a pair
$(h_0,j_0)$ such that $m<\omega_{h_0,j_0}<h_0\}$. Hence, for any
$\omega\in B_{m,n}$, the pair $(h_0,j_0)$ satisfies the inequality
$h_0>m$. Since the columns of $\omega$ are non-increasing, for all
$k\le m$, we have $\omega_{k,j_0}\ge\omega_{h_0,j_0}$ and
$\omega_{k,j_0}\ge k$, whence at least $m$ parts of $\omega$ are
greater than $m$. So, we observe that $A_{m,n}\subset
B_{m,n}\subset\{X_{m,n}>m\}$ and thus
\begin{equation}\label{tail}
\mathbb{P}(A_{m,n})\le\mathbb{P}(X_{m,n}>m).
\end{equation}
By Theorem \ref{t2}, as $n\to\infty$, the upper bound in
(\ref{tail}) deals with a right tail of a Poisson distribution. In
addition, note that by (\ref{mp}) we have $m\to\infty$ as
$n\to\infty$. Therefore we conclude that, under the assumptions of
Theorem \ref{t2}, the right hand side of (\ref{tail}) tends to
$0$. Hence $\mathbb{P}(A_{m,n})\to 0$, which means that
$Z_{m,n}-X_{m,n}\to 0$ in probability. Now, from the
representation $Z_{m,n}=(Z_{m,n}-X_{m,n})+X_{m,n}$ it follows that
$Z_{m,n}$ and $X_{m,n}$ have one and the same limiting
distribution as $n\to\infty$, which completes the proof of the
corollary. $\Box$

\section{Concluding Remarks}
\setcounter{equation}{0}

In this work, we have proposed an approach for analysis of
statistics related to counts of large and small parts in a random
and uniformly chosen plane partition. The concept of conjugate
trace of a plane partition, introduced by Stanley \cite{S73},
plays an important role in our study. We have decomposed the
conjugate trace into the sum of the counts of parts $1,2,\dots$
and then we have separated these counts by the formal variables
$y_1,y_2,\dots$ in the underlying generating function. Cutting
this sum until the $m$th row of the plane partition array shadowed
by the conjugate trace, for suitable values of $m$ (see
(\ref{mp})), we are able to study the asymptotic behavior of
several particular statistics of counts of large part sizes. As an
illustration, we have proved a Poisson approximation for a
statistic, which has a natural analogue in the case of linear
integer partitions. Finally, we have shown that the restriction
given by the conjugate trace of a plane partition of $n$ can be
removed since it is bounded by a tail of a Poisson distribution
and tends to $0$ as $n\to\infty$. We believe that our approach
could be also applied to further studies in this direction (for
instance, to establish central and local limit theorems for other
similar plane partition statistics).

Our last remark is related to the one-dimensional case of linear
integer partitions. Fristedt's method of study of integer
partition statistics \cite{F93} is purely probabilistic. It
transfers the joint probability distribution of the part counts of
a random and uniformly chosen one-dimensional partition into the
joint conditional distribution of independent and geometrically
distributed random variables. We are able to give an alternative
proof of Theorem \ref{t1} based on generating function identities.
Below, for the sake of completeness, we will briefly sketch this
proof.

First, using general results from \cite[Chapter 1]{A76}, it is not
difficult to show that
$$
1+\sum_{n=1}^\infty p(n)\mathcal{E}(y^{\xi_{m,n}})x^n
=P(x)\varphi_m(x,y),
$$
where
$$
\varphi_m(x,y) =\prod_{j>m}\frac{1-x^j}{1-yx^j}.
$$
We recall that $\mathcal{E}$ denotes the expectation with respect
to the uniform probability measure $\mathcal{P}$ on the set
$\Lambda(n)$ of linear integer partitions of $n$, and that $P(x)$,
the generating function of the sequence $\{p(n)\}_{n\ge 1}$,
satisfies (\ref{pgf}). Then, one can apply the asymptotic scheme
for analysis of linear integer partition statistics, proposed by
Grabner et al. \cite[Theorem 2.2]{GKW14} and based on the
classical saddle point method, to show that, for
$$
d_n^\prime =\frac{\pi}{\sqrt{6n}} -\frac{1}{4n} +O(n^{-1-\alpha}),
\quad \alpha>0,
$$ and $m$ defined by (\ref{mpf}), we have
$$
\lim_{n\to\infty} \varphi_m(e^{-d_n^\prime},y)
=\exp{(e^{-c}(y-1))}
$$
uniformly for $y\in [0,1)$. The last exponent is the probability
generating function of the Poisson distribution with expectation
$e^{-c}$, which completes the proof of Theorem \ref{t1}.

\section*{Acknowledgements}

The author is grateful to the anonymous referees and the editor
for their valuable comments and especially for suggesting a number
of clarifications and corrections on an earlier draft of this
paper.

\end{document}